\documentclass[a4paper]{article}
\usepackage{cite}
\usepackage[all]{xy}
\usepackage[latin1]{inputenc}      
\usepackage[dvips]{graphics,graphicx}
\usepackage{amsfonts,amssymb,amsmath,color,mathrsfs, amstext}
\usepackage{amsbsy, amsopn, amscd, amsxtra, amsthm,authblk}
\usepackage{enumerate,algorithmicx,algorithm}
\usepackage{algpseudocode}
\usepackage{upref}
\usepackage{geometry,wrapfig}
\usepackage{amsthm,amsmath,amssymb}
\usepackage{relsize}
\usepackage{mathrsfs}
\usepackage{bm}
\usepackage{ulem}
\usepackage{exscale}
\usepackage{tabularx}
\usepackage{threeparttable,multirow,dcolumn,booktabs}
\usepackage{algorithmicx,algorithm}
\usepackage{tabularx}
\usepackage{caption}
\usepackage[displaymath]{lineno}

\usepackage{float}
\usepackage{bm}
\usepackage{caption}
\usepackage{yhmath}
\usepackage{booktabs}
\usepackage{subcaption}
\usepackage{multirow}
\usepackage{makecell}
\usepackage[colorlinks,
            linkcolor=red,
            anchorcolor=red,
            citecolor=red
            ]{hyperref}

\numberwithin{equation}{section}

\numberwithin{equation}{section}

\newcommand{\bs}[1]{\boldsymbol{#1}}
\newcommand{\bmu}{\bs{\mu}}

\newcommand{\bk}{\bs{k}}
\newcommand{\calN}{\mathcal{N}}
\newcommand{\calA}{\mathcal{A}}

\numberwithin{equation}{section}
\usepackage{epstopdf}
\usepackage{pdfpages}

\graphicspath{{./Figs/}}

\begin{document}

\title{MCMS-RBM: Multi-Component Multi-State Reduced Basis Method toward Efficient Transition Pathway Identification for Crystals and Quasicrystals}

\author{Yajie Ji\thanks{School of Mathematical Sciences, Shanghai Jiao Tong University, Shanghai 200240, China. Email: {\tt jiyajie595@sjtu.edu.cn}. },\quad
{Lijie Ji}\thanks{School of Mathematical Sciences, Shanghai Jiao Tong University, Shanghai 200240, China. Corresponding author. Email: {\tt sjtujidreamer@sjtu.edu.cn}. L. Ji acknowledges the support from NSFC grant No. 12201403 and China Postdoctoral Science Foundation No. 2021M702141.},\quad Yanlai Chen\footnote{Department of Mathematics, University of Massachusetts Dartmouth, 285 Old Westport Road, North Dartmouth, MA 02747, USA. Email: {\tt{yanlai.chen@umassd.edu}}. Y. Chen is partially supported by National Science Foundation grant DMS-2208277 and by the UMass Dartmouth MUST program, N00014-20-1-2849, established by Dr. Ramprasad Balasubramanian via sponsorship from the Office of Naval Research.}, \quad  
Zhenli Xu\footnote{School of Mathematical Sciences, CMA-Shanghai, and MOE-LSC, Shanghai Jiao Tong University, Shanghai 200240, China. Email: {\tt xuzl@sjtu.edu.cn}. Z. Xu acknowledges the support
from the NSFC (grant No. 12071288). The authors also acknowledge the support from the HPC center of Shanghai Jiao Tong University.}}

\date{}
\makeatother

\maketitle

\begin{abstract}
Due to quasicrystals having long-range orientational order but without translational symmetry, traditional numerical methods usually suffer when applied as is. In the past decade, the projection method has emerged as a prominent solver for quasiperiodic problems. Transforming them into a higher dimensional but periodic ones, the projection method facilitates the application of the fast Fourier transform. However, the computational complexity inevitably becomes high which  significantly impedes e.g. the generation of the phase diagram since a high-fidelity simulation of a problem whose dimension is doubled must be performed for numerous times.

To address the computational challenge of quasiperiodic problems based on the projection method, this paper proposes a multi-component multi-state reduced basis method (MCMS-RBM). Featuring multiple components with each providing reduction functionality for one branch of the problem induced by one part of the parameter domain, the MCMS-RBM does not resort to the parameter domain configurations (e.g. phase diagrams) \textit{a priori}. It enriches each component in a greedy fashion via a phase-transition guided exploration of the multiple states inherent to the problem. Adopting the empirical interpolation method, the resulting online-efficient method  vastly accelerates the generation of a delicate phase diagram to a matter of minutes for a parametrized two-turn-four dimensional Lifshitz-Petrich model with two length scales. Moreover, it furnishes surrogate and equally accurate field variables anywhere in the parameter domain.


\end{abstract}

\section{Introduction}

In the 1980s, Shechtman \textit{et al.} \cite{shechtman1984metallic} observed a metallic phase with long-range orientational order in a rapidly cooling Al-Mn alloy.  
  Unlike the periodic structures that feature translational symmetries or 1-, 2-, 3-, 4-, and 6-fold rotational symmetries,  
  this new structure has a 5-fold rotational symmetries. Later, researchers coined this new long-range ordered structure ``quasicrystals" \cite{levine1984quasicrystals}.
Since the discovery of the 5-fold quasicrystals, more different structures with 5-, 6-, 8-, 12-, and 20-fold symmetries have emerged in various metallic alloys \cite{tsai2008icosahedral,steurer2004twenty}. There have also been certain soft quasicrystals in soft matter systems \cite{percec2009self,percec2009self2,
hayashida2007polymeric,talapin2009quasicrystalline,zeng2004supramolecular,zhang2012dodecagonal}. Since their discoveries, these quasicrystals are widely used in materials science, thermal engineering, metallurgical engineering, photonics, and energy storage\cite{tsai2008icosahedral,zeng2004supramolecular,zhang2012dodecagonal}.

\begin{figure}
\centering
\includegraphics[scale=0.58]{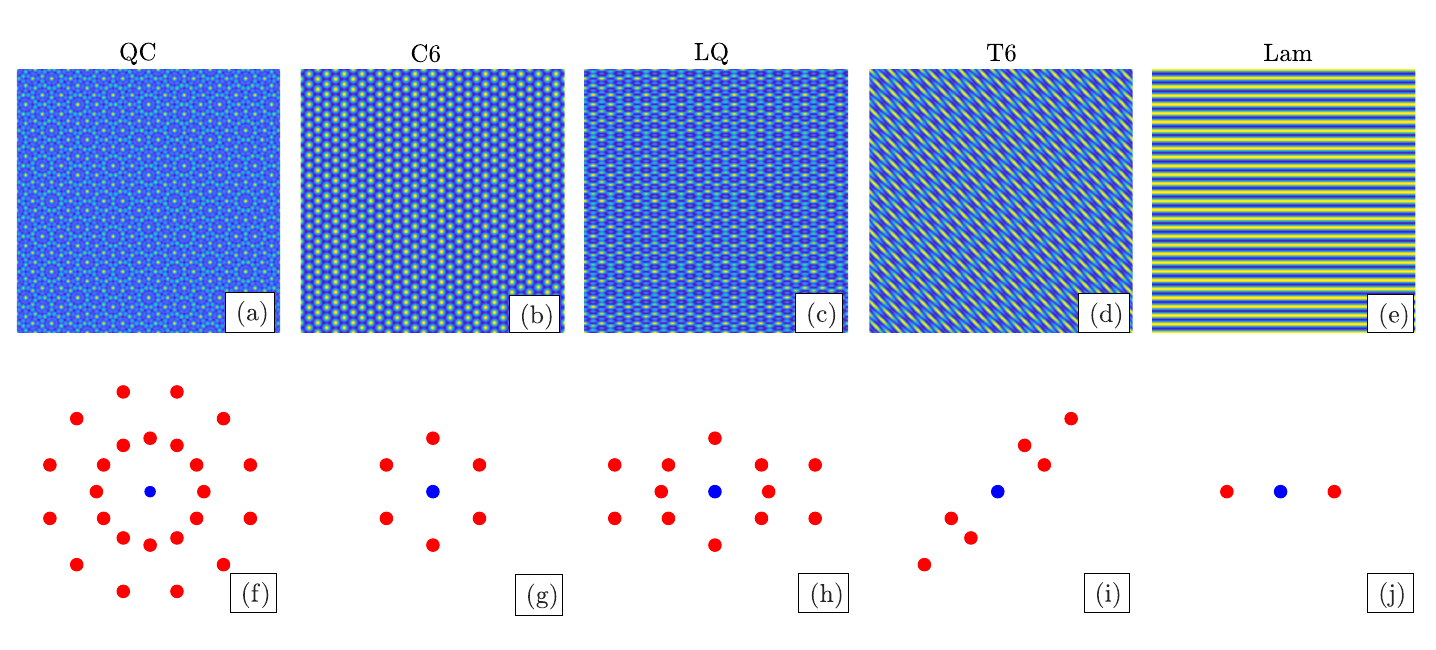}
\caption{Order parameters (Top) and the corresponding prominent diffraction patterns in the reciprocal space (Bottom) for QC (af), C6(bg), LQ(ch), T6(di), and Lam(ej) states. }
\label{PNAS_ref}
\end{figure}

In order to understand the formation, stability and various physical properties of different quasicrystals, some theoretical and numerical methods have been proposed. For the theoretical model, 
Lifshitz and Petrich first proposed the Lifshitz-Petrich (LP) model to explain the 12-fold symmetry excited by dual-frequency filtering in the Faraday experiment\cite{lifshitz1997theoretical}. Subsequently, the LP model have been used to study the stability of the two-dimensional 5-, 8-, and 10-fold 
quasicrystals with two characteristic wavelength scales \cite{jiang2015stability}. The LP model represents a coarse-grained mean field theory. It assumes that the free energy of the system can be represented as a function of the order parameters. There are two such parameters, with one representing the temperature and the other delineating the asymmetry of the order parameter. When these two parameters vary, the quasicrystals exhibit a rich phase behavior containing a number of equilibrium ordered phases. This phase diagram, when captured well, can be used to study the transition path between different structures. For example, the LP model describing the 12-fold quasicrystals in two dimensions also includes the 6-fold crystalline state (C6), the lamellar quasicrystalline state (LQ), the transformed 6-fold crystalline state (T6) and the Lamella state (Lam) \cite{yin2021transition}. The structures of these stable states are shown in real and the so-called reciprocal spaces, respectively in Fig. \ref{PNAS_ref}. 

It is therefore imperative to simulate these quasicrystals accurately across the entire parameter domain. That proves extremely challenging for two reasons. 
First, each simulation is delicate and costly. Unlike periodic structures that are translation invariant and therefore can be calculated in a unit cell with periodic boundary conditions, quasicrystals are rotationally invariant but not translationally invariant.  This makes it difficult to determine the computational domain and boundary conditions.
Second, the need of repeated simulations exacerbate the situation. To resolve the phase diagram even on a relatively small parameter domain, thousands of simulations are needed.

There are two popular methods to overcome the first challenge:
crystaline approximation method (CAM) \cite{lifshitz1997theoretical,zhang2008efficient}, and the projection method (PM) \cite {jiang2014numerical}. The CAM approximates the quasicrystals with a periodic structure, and the size of the computational domain increases rapidly as the accuracy of approximation becomes higher. This has been systematically illustrated in several papers  and readers can refer to \cite{jiang2014numerical}.
The PM utilizes the fact that the reciprocal vectors of a quasicrystal in a lower dimensional space can be approximated by a linear combination of basic reciprocal vectors in a higher dimensional space. This technique renders the quasicrystal periodic in a higher-dimensional space \cite{walter2009crystallography}.
With this idea, the Fourier expansion approach can be employed for the quasiperiodic systems \cite{jiang2014numerical,jiang2015stability,
jiang2017stability,zhou2019plane,gao2023pythagoras}.
This projection method is performed in high dimensional space and can use the fast Fourier transform. 
For the two dimensional LP model with two length scales, the phase diagram is quite complicated due to the wide range of parameter values and large number of possible stable states. It is severely time consuming if one wants to generate the phase diagram for a wide range of parameters 
accurately \cite{Quarteroni2015}. This is an especially onerous task if the physical and/or parametric domain are of high dimensions.  Although adaptive method exists \cite{Softwarejiang} that can generate the phase diagram without having to resolve the full parameter domain, it fails to produce the field variables which are needed e.g. in controlling the self-assembly of quasicrystals and a variety of other desired structures in practical experimental realizations.\cite{provatas2011phase,chen2002phase,
barkan2011stability,lifshitz2007soft}. 
Moreover, if there are five possible stable states for each unknown parameter, one needs to solve the LP model five times with respect to five different initial values, and then choose the one having the minimum free energy as the stable state. Therefore, the pursuit of efficient numerical algorithms for the parametric LP model has emerged as a prominent research focus.

To achieve that goal while dealing with the second challenge, we propose in this paper a multi-component multi-state reduced basis method (MCMS-RBM) as a generic framework for reduced order modeling for parametric problems whose solution has multiple states across the parameter domain. The RBM \cite{Quarteroni2015,rozza2008reduced,
hesthaven2016certified} is a projection-based model order reduction technique that provides a mathematically accurate surrogate solution in a highly efficient manner, and  capable of reducing the computational complexity of the full order model (FOM) by several orders of magnitude after an offline learning stage. It was first introduced for nonlinear structure problem in 1970s \cite{noor1980reduced} and has been later analyzed and extended to solve many problems such as linear evolutionary equation \cite{haasdonk2008reduced}, viscous Burgers equation\cite{veroy2003reduced}, the Naiver-Stokes equations \cite{deparis2009reduced} among others.  Its extension to the Reduced Over Collocation setting \cite{ChenJiNarayanXu2020, ChenSigalJiMaday2021, chen2022hyper} makes available a robust and efficient implementation for the nonlinear and non-affine setting.   
The offline-online decomposition, often assisted by the empirical interpolation method (EIM) \cite{barrault2004empirical,grepl2007efficient}, is a critical approach to realize online efficiency meaning the online solver is independent of the degrees of freedom of the FOM. 
The proposed MCMS-RBM has two prominent features in comparison to the standard RBM. First, it has multiple components with each providing reduction functionality for one branch of problem induced by one part of the parameter domain. Second, without resorting to the parameter domain configurations (e.g. phase diagrams) \textit{a priori}, it enriches each component in a greedy fashion. With the existence of multiple stable states and the occurrences of phase transitions, it is difficult to construct precise low dimesnional RB spaces for each candidate state. The MCMS-RBM overcomes this difficulty by leveraging the structures of the prominent reciprocal vectors of all possible stable states and designing a phase transition indicator. This indicator guides the greedy algorithm to explore the multiple states inherent to the quasicrytal problem.

The rest of this paper is organized as follows. In Section \ref{sec:LP_PM}, we introduce the parametrized LP model, the projection method, pseudospectral method used to obtain the high-fidelity solution, the phase diagram, and the criteria of phase transition. The key components of the MCMS-RBM, including the online and offline procedures and the EIM process, are introduced in Section \ref{sec:mcms}. We then present numerical results in Section \ref{sec:numerics} to demonstrate the efficiency and accuracy of the proposed MCMS-RBM. Finally, concluding remarks are drawn in Section \ref{sec:conclusion}.

\section{Lifshitz-Petrich model and the projection method}
\label{sec:LP_PM}
In this section, we introduce the  Lifshitz-Petrich model with two characteristic length scales and the phase-steady full order solutions based on the projection method.

\subsection{Lifshitz-Petrich model}

The LP model is based on the Swift-Hohenberg model \cite{swift1977hydrodynamic} and the Landau-Brazovskii model \cite{brazovskii1975phase}, which are widely used in the study of materials science and polymeric systems, respectively. The LP model extends one wavelength scale of the Swift-Hohenberg equation to two characteristic wavelength scales. 
The scalar order parameter $\phi(\bm{r})$ describes how perfectly the molecules are aligned. It minimizes 
the corresponding free energy functional which is defined as
\begin{equation}
\mathcal{F}(\phi; c, q, \bmu)=\int_V d r\left\{\frac{c}{2}\left|\left(\nabla^2+1^2\right)\left(\nabla^2+q^2\right) \phi\right|^2-\frac{\varepsilon}{2} \phi^2-\frac{\alpha}{3} \phi^3+\frac{1}{4} \phi^4\right\}, \label{phi_energy}
\end{equation}
where $\bm r \in \mathcal{R}^d$ with $d=2$, $V$ is the system volume, $c$ is an energy penalty parameter to ensure that the principle reciprocal vectors of structures is located on $|\bm{k}|=1$ and $|\bm{k}|=q$, with $q$ being an irrational number depending on the symmetry, $\varepsilon$ is the reduced temperature and $\alpha$ is a phenomenological parameter. For simplicity, we define $\bmu:=[\varepsilon,\alpha]$, a two dimensional parameter vector.

For a given parameter $\bmu$, the candidate stable states are the local minima of the free energy functional, that is, the solutions of the 
Euler-Language equation
\begin{equation}
\frac{\delta \mathcal{F}}{\delta \phi}= 0.
\end{equation}
This is a eighth-order nonlinear partial differential equation. To solve it, one can use gradient flow method \cite{jiang2020efficient,jiang2015stability}
\begin{equation}
\frac{\partial \phi}{\partial t}=-\frac{\delta \mathcal{F}}{\delta \phi}=-\left[c\left(\nabla^2+1\right)^2\left(\nabla^2+q^2\right)^2 \phi-\varepsilon \phi-\alpha \phi^2+\phi^3\right],
\end{equation}
which is then discretized by the following implicit-explicit scheme
\begin{equation}
\frac{\phi^{n+1}-\phi^n}{\Delta t}=-\left[c\left(\nabla^2+1\right)^2\left(\nabla^2+q^2\right)^2 \phi^{n+1}-\varepsilon \phi^n-\alpha\left(\phi^n\right)^2+\left(\phi^n\right)^3\right].
\label{eq:semi_implicit}
\end{equation}

Depending on the values of $\bmu$, solutions of this eighth-order nonlinear partial differential equation lead to quasicrystals or periodic structures. For the latter, there are many fast algorithms. Specifically, the pseudospectral method achieves efficiency by evaluating the gradient terms in the Fourier space and the nonlinear terms in the physical space. On the other hand, the quasiperiodic structure cannot be solved directly by these classical methods suitable for the periodic structure. We adopt the projection method, the topic of the next section.

\subsection{Phase-steady full order solutions based on the projection method}
\label{sec:PM}
In the PM \cite{jiang2014numerical}, a quasicrystal is first computed in a (higher-dimensional) reciprocal space as a periodic structure which is then projected back to the lower-dimensional space through the projection matrix. We provide a brief review of this approach. 

One can represent the reciprocal vectors $\bm{k} \in \mathbb{R}^d$ of a $d$-dimensional quasicrystal as \cite{chaikin1995principles}
\begin{align*}
\bm{k}=h_1 \bm{p}_1^*+h_2 \bm{p}_2^*+\cdots+h_n \bm{p}_n^*, \quad h_i \in \mathbb{Z},
\end{align*}
with $\bm{p}_i^* \in \mathbb{R}^d$ having $\mathbb{Z}$-rank of $n$\footnote{The only $\{h_i \in \mathbb{Z}\}$ leading to $h_1 \bm{p}_1^*+h_2 \bm{p}_2^*+\cdots+h_n \bm{p}_n^* = 0$ is $h_i \equiv 0$.}  ($n>d$). Different choices of the coefficient vector
\[
\bm{h} \triangleq \{h_1, \cdots, h_n\}
\]
(e.g. setting some of them to be zero and enforcing constraints on the others) lead to different quasicrystal patterns, see Fig. \ref{PNAS_ref}. The PM finds proper $n$-dimensional vectors $\bm{b}_i^*, 1\leq i\leq n,$ being the primitive reciprocal vectors of the $n$-dimensional reciprocal space. The reciprocal vector of $n$-dimensional periodic structure is then 
$$\bm{H}=h_1 \bm{b}_1^*+h_2 \bm{b}_2^*+\cdots+h_n \bm{b}_n^*,\quad h_i \in \mathbb{Z}.$$
We denote by $\mathcal{S} \in \mathbb{R}^{d \times n}$ a projection matrix satisfying $\bm{p}_i^* = \mathcal{S}\bm{b}_i^*$.

\begin{figure}
\centering
\includegraphics[scale=0.42]{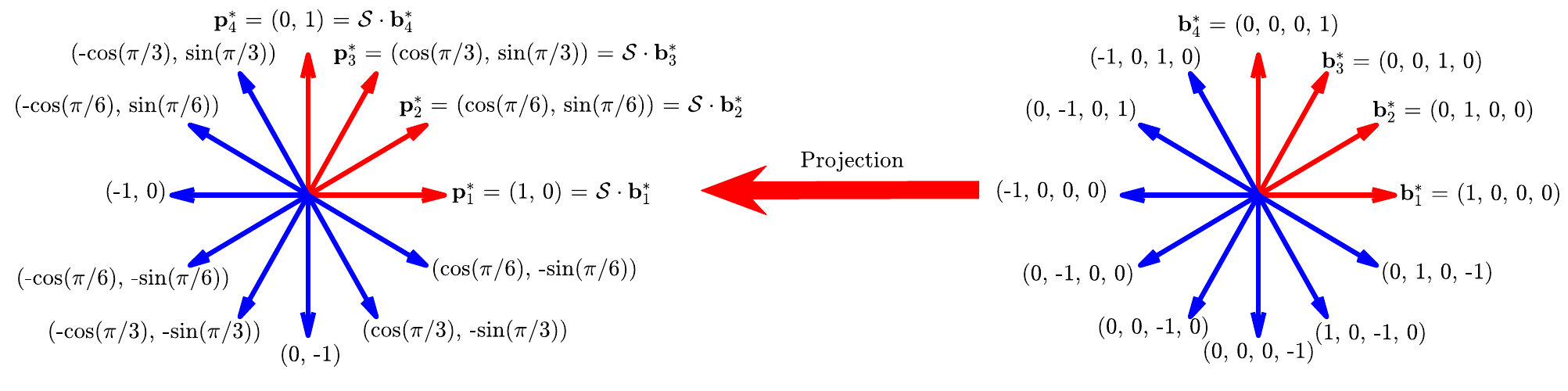}
\caption{Reciprocal lattice vectors  of 12-fold rotational symmetry in 2- and 4-dimensional space.  The two sets of vectors are connected through the projection matrix $\mathcal{S}$. }
\label{fig:reciprocalbasis}
\end{figure}
This means that the $d$-dimensional quasicrystal is a periodic structure in the $n$-dimensional space. In this paper, we focus on the case that $n=4$ and $d=2$ and the  quasicrystal is 12-fold (i.e. we take $q=2\cos (\pi/12)$ in the free energy functional $\mathcal{F}(\phi(\bm{r}); c, q, \bmu)$ \eqref{phi_energy}.   
The primitive reciprocal vectors are $\bm{p}_1^*=(1,0)$ and $\bm{p}_4^*=(0,1)$, see Fig. \ref{fig:reciprocalbasis}. 
However, some reciprocal 
vectors of the 12-fold case cannot be represented as a linear combination of $\bm{p}_1^*$ and $\bm{p}_4^*$ with integral coefficients. We therefore adopt 
$$\bm{p}_1^*=(1,0),\quad \bm{p}_2^*=(\cos(\pi/6),\sin(\pi/6)),\quad \bm{p}_3^*=(\cos(\pi/3),\sin(\pi/3)),\quad \bm{p}_4^*=(0,1)$$ of $\mathbb{Z}$-rank 4, and the projection matrix 
$$\mathcal{S}=\left(\begin{array}{llll}
1 & \cos (\pi / 6) & \cos (\pi / 3) & 0 \\
0 & \sin (\pi / 6) & \sin (\pi / 3) & 1
\end{array}\right).$$

 With this idea, the Fourier expansion for the $d$-dimendional quasiperiodic function is given by
\begin{equation}
\phi(\bm{r})=\sum_{\bm{H}} \widehat{\phi}(\bm{H}) e^{i\left[(\mathcal{S}\cdot \bm{H})^T \cdot \bm{r}\right]},\quad \bm{r} \in \mathbb{R}^d,\quad \bm{H} \in \mathbb{R}^n.
\label{PM_solver}
\end{equation}
Denoting by $g_k^T$ the $k^{\rm th}$ row of $\mathcal{S} \cdot \bm{H}$, we have that 
$$\mathcal{S} \cdot \bm{H}=\left(\sum_{i=1}^n s_{1 i} \sum_{j=1}^n h_j b_{j i}^*, \cdots, \sum_{i=1}^n s_{d i} \sum_{j=1}^n h_j b_{j i}^*\right)^T \triangleq\left(g_1, \cdots, g_d\right)^T,\quad h_j \in \mathbb{Z},$$
with $b_{j,i}^*, j=1, \cdots, n$ being the components of $b_i^*$.
The LP free energy functional then becomes
\begin{equation}
\begin{aligned}
\mathcal{F}(\phi(\bm{r}); c, q, \bmu)= & \frac{1}{2} \sum_{\bm{H}_1+\bm{H}_2=0}\left\{c(1-\sum_{k=1}^d g_k^2)^2(q^2-\sum_{k=1}^d g_k^2)^2-\varepsilon\right\} \widehat{\phi}\left(\bm{H}_1\right) \widehat{\phi}\left(\bm{H}_2\right) \\
& -\frac{\alpha}{3} \sum_{\bm{H}_1+\bm{H}_2+\bm{H}_3=0} \widehat{\phi}\left(\bm{H}_1\right) \widehat{\phi}\left(\bm{H}_2\right) \widehat{\phi}\left(\bm{H}_3\right)\\
&+\frac{1}{4} \sum_{\bm{H}_1+\bm{H}_2+\bm{H}_3+\bm{H}_4=0} \widehat{\phi}\left(\bm{H}_1\right) \widehat{\phi}\left(\bm{H}_2\right) \widehat{\phi}\left(\bm{H}_3\right) \widehat{\phi}\left(\bm{H}_4\right),
\end{aligned}
\label{eq:high_free_energy}
\end{equation}
Substituting Eq. \eqref{PM_solver} into Eq. \eqref{eq:semi_implicit} and using Eq. \eqref{eq:high_free_energy}, one obtains
\begin{equation}
\begin{aligned}
&\left(\frac{1}{\Delta t}+c\left(1-\sum_{k=1}^d g_k^2\right)^2\left(q^2-\sum_{k=1}^d g_k^2\right)^2\right) \widehat{\phi}_{t+\Delta t}(\bm{H}) =\left(\frac{1}{\Delta t}+\varepsilon\right) \widehat{\phi}_t(\mathbf{H})\\
&+\alpha\sum_{\bm{H}_1+\bm{H}_2=\bm{H}} \widehat{\phi}_t\left(\bm{H}_1\right) \widehat{\phi}_t\left(\bm{H}_2\right)-\sum_{\bm{H}_1+\bm{H}_2+\bm{H}_3=\bm{H}} \widehat{\phi}_t\left(\bm{H}_1\right) \widehat{\phi}_t\left(\bm{H}_2\right) \widehat{\phi}_t\left(\bm{H}_3\right),
\end{aligned}
\label{eq:semi_implicit_pm}
\end{equation}
where $\Delta t$ is the temporal step, $\widehat{\phi}_{t+\Delta t}$ and $\widehat{\phi}_{t}$ represent the Fourier coefficients at time $t+\Delta t$ and $t$, respectively.
A direct evaluation of the convolution terms of \eqref{eq:semi_implicit_pm} are expensive. Instead, one can calculate these nonlinear terms in the physical space and then perform FFT to derive the corresponding Fourier coefficients. Therefore,
the computational complexity of the PM is $$O\left( N_{t}\cdot \calN  \log \calN \right),$$
where $N_{t}$ is the number of time iterations, and $\calN= \left(N_{\bm{H}}\right)^n$ with $N_{\bm{H}}$ being the degrees of freedom of pseudospectral method in each dimension.

\subsection{Phase diagram, phase transition and multiple phase-steady solutions}

\label{sec:pdpt}

The phase diagram is a quantitative and graphical representation of the stability and interconversion relationships of various metastable/stable phases of a material under different conditions, e.g. with different temperature $\varepsilon$ and phenomenological parameter $\alpha$. The phase field model can be used to not only simulate the phase transformation and microstructure changes during the processing and handling of materials, but also predict the possible emergence of new materials or novel phases. However, it is quite time consuming to generate the phase diagram for a wide range of parameter values.

For each value of parameter $\bmu:=[\varepsilon,\alpha]$, due to the existence of the multiple stable solutions corresponding to the different prominent reciprocal vectors, one needs to solve Eq. \eqref{eq:semi_implicit_pm} five times with five different initial values $\widehat{\phi}_0(\bm{\mathcal H}_{\rm S})$ corresponding to five candidate states $\bm{\mathcal H}_{\rm S}$ for $\rm S \in \{{\rm QC, C6, LQ, T6, Lam}\}$. 
The iteration initialized specifically based on the reciprocal vectors of each state allows for a rapid convergence of the gradient flow equation. 
Indeed, we denote by 
$$\widehat{\phi}(\cdot; \bm{\mathcal H}_{\rm S}, \bmu)$$ 
the steady-state solution for $\bmu$ with the initial value given by $\widehat{\phi}_0(\bm{\mathcal H}_{\rm S})$. As indicated in Figure \ref{PNAS_ref}(f), there are 24 prominent reciprocal vectors for $\rm S = QC$. These 2D reciprocal vectors can be transformed into four dimensional space and they are related by the projection matrix $\mathcal{S}$, see Figure \ref{fig:reciprocalbasis}. For a more intuitive illustration, we display all the 24 prominent reciprocal vectors in four-dimensional space in Table \ref{table:initialindex}, and the remaining sets $\bm{\mathcal H}_{\rm S}$'s for $\rm S \in \{C6, LQ, T6, Lam\}$ are defined as follows. For $\rm S = C6$, $\bm{\mathcal H}_{\rm S}$ contains the 6 reciprocal vectors displayed in bold.  For $\rm S = LQ$, it contains the 12 reciprocal vectors underlined. For $\rm S = T6$, the 6 prominent reciprocal vectors are displayed with dash lines. The 2 reciprocal vectors for $\rm S = Lam$ are displayed with wavy lines. For a rapid convergence, the Fourier coefficients of these five sets of reciprocal vectors are initialized with nonzero values
\begin{equation}
\widehat{\phi}_0(\bm{\mathcal H}_{\rm S})=u_0,
\end{equation}
where $u_0$ is a given constant.
\begin{table}[thbp]
\renewcommand\arraystretch{2.0}
	\centering
			\begin{tabular}{c|cccccc}
		\hline
\multirow{2}{*}{$|\mathcal{S} \cdot \mathbf{H}| =1$} 
&\uwave{\underline{(0 1 0 0)}}& \underline{\textbf{\dashuline{(0 0 1 0)}}}& \dashuline{(0 0 0 1)} &\underline{\textbf{(-1 0 1 0)}}& (0 -1 0 1)& \underline{\textbf{(-1 0 0 0)}} \\ 
~~ &\uwave{\underline{(0 -1 0 0)}} &\underline{\textbf{\dashuline{(0 0 -1 0)}}}& \dashuline{(0 0 0 -1)}  &\underline{\textbf{(1 0 -1 0)} }&(0 1 0 -1) &\underline{\textbf{(1 0 0 0)}}\\ \hline
\multirow{2}{*}{$|\mathcal{S} \cdot \mathbf{H}| =q$} 
&\underline{(1 1 0 0)}& \dashuline{(0 1 1 0)}& (0 0 1 1)&(-1 0 1 1)& (-1 -1 1 1)& \underline{(-1 -1 0 1)} \\ 
~~ &\underline{(-1 -1 0 0)} &\dashuline{(0 -1 -1 0)}& (0 0 -1 -1)  &(1 0 -1 -1) &(1 1 -1 -1) &\underline{(1 1 0 -1)}\\  \hline
	\end{tabular}
	\caption{Prominent reciprocal vectors in four dimensional space with nonzero initial Fourier coefficients.} \label{table:initialindex}
\end{table}

The set of \textit{multiple phase-steady solutions} (PSS) corresponding to $\bmu$ is then 
\begin{equation}
\label{eq:mu_pss}
\widehat{\Phi}(\cdot; \bmu) \triangleq \bigcup_{{\rm S} = {\rm QC}}^{\rm Lam}\left\{\widehat{\phi}(\cdot; \bm{\mathcal H}_{\rm S}, \bmu)\right\}
\end{equation}
Here, for simplicity, we define $\left\{\widehat{\phi}(\cdot; \bm{\mathcal H}_{\rm S}, \bmu)\right\}$ to be the empty set when $\widehat{\phi}$ goes through a phase transition with the given initial value $\widehat{\phi}_0(\bm{\mathcal H}_{\rm S})$. 
The rationale is that since each parameter leads to a stable state solution without phase transitions, the solutions that undergo phase transitions during the evolution process can be readily discarded.

The existence of multiple convergent solutions for the same $\bmu$ serves two purposes. 
On one hand for the single-query setting, the state within $\widehat{\Phi}(\cdot; \bmu)$ that leads to the smallest free energy functional is the stable state for the queried parameter. 
On the other hand for the multi-query setting, the construction of the multiple components of our proposed MCMS-RBM takes advantage of the existence of the multiple solutions corresponding to the multiple states. The many-to-many pattern between components and states, a main novelty of the MCMS-RBM, enables the quick and simultaneous enriching of the reduced spaces with limited FOM queries. The detailed high fidelity solver for 
\[
\bmu \mapsto \widehat{\Phi}(\cdot; \bmu)
\]
is shown in Algorithm \ref{alg:Algorithm:PMsolver}. 
\begin{algorithm}[h]
  \caption{FOM for the LP model, $\bmu \mapsto \widehat{\Phi}(\cdot; \bmu)$}
  \label{alg:Algorithm:PMsolver}  
  \begin{algorithmic}[1]  
    \vspace{0.5ex}
\State \textbf{Input:} parameter $\bmu$, tolerance $\operatorname{tol}$, temporal step $\Delta t$, and maximum iteration $T$
\State Set $\widehat{\Phi}(\cdot; \bmu) = \{\}$
\State \mbox{\textbf{For}} ${\rm S} \in \{{\rm QC, C6, LQ, T6, Lam}\}$
    \State \quad Initialize the Fourier coefficients $\widehat{\phi}_0$ corresponding to $\bm{\mathcal H}_{\rm S}$
    \State \quad Initialize $\operatorname{res} = 1$, ${\rm PTI}_{S}=0$, and $t=0$. Calculate $\phi_0 = \mathscr{F}^{-1}(\widehat{\phi}_0)$
    \State \quad \mbox{\textbf{While}} {$\operatorname{res}\geq \operatorname{tol}$ \& $ t\leq T$ \& ${\rm PTI}_{\rm S}=0$}   
    \State \quad\quad Obtain $\widehat{\phi}_{t+\Delta t}$ by Eq. \eqref{eq:semi_implicit_pm} using $\widehat{\phi}_{t}$ and ${\phi}_{t}$
    \State \quad\quad $\operatorname{res}=||\widehat{\phi}_{t+\Delta t}-\widehat{\phi}_t||_2$
    \State \quad\quad $t = t+ \Delta t$, $\phi_t = \mathscr{F}^{-1}(\widehat{\phi}_t)$
    \State \quad \quad ${\rm PTI}_{\rm S}={\rm PTI}(\widehat{\phi}_0, \widehat{\phi}_t)$
    \State \quad\mbox{\textbf{End While}}
    \State \quad If ${\rm PTI}_{\rm S}=0$, set $\widehat{\Phi}(\cdot; \bmu) = \widehat{\Phi}(\cdot; \bmu) \cup \{\widehat{\phi}_t\}$
\State \mbox{\textbf{End For}}
    \State \textbf{Output:} High fidelity PSS set $\widehat{\Phi}(\cdot; \bmu)$
  \end{algorithmic}  
\end{algorithm}

This algorithm utilizes a phase transition indicator (PTI) which is given in Algorithm  \ref{alg:Algorithm:phasetransit}. It is inspired by the observation that, if the initial state is not in the stable state (out of the $5$ possible states) corresponding to the given parameter, the evolution may or may not undergo a transition to other states. When the transition happens, it deteriorates the low-rank nature of the corresponding branch of the solution manifold. We therefore discard the convergent solution whenever phase transition occurs in the evolution process. 
The detection of such transitions is made possible by the realization that solutions sharing the same structure maintain consistent characteristics in their corresponding \textit{spectral signature} (the set of spectral coefficients whose magnitudes are above a certain tolerance) throughout the entire evolution. The ``emergence'' or ``disappearance'' of a spectral mode therefore signifies a phase transition. 
\begin{algorithm}[h]
  \caption{Phase Transition Indicator, 
  ${\rm PTI}~ (\widehat{\phi}_0, \widehat{\phi}_t)$}
  \label{alg:Algorithm:phasetransit}  
  \begin{algorithmic}[1]  
    \vspace{0.5ex}
\State {\bf Input:} Reference state $\widehat{\phi}_0$, candidate $\widehat{\phi}_t$, and tolerance $\delta$
\State Calculate the locations of the elements of $\widehat{\phi}_0$ with $|\widehat{\phi}_0|>\delta$
\State Calculate the locations of the elements of $\widehat{\phi}_t$ with $|\widehat{\phi}_t|>\delta$
\State If they coincide, then ${\rm PTI} =0$, otherwise, ${\rm PTI} = 1$
  \end{algorithmic}  
\end{algorithm}

\section{The MCMS-RBM}

\label{sec:mcms}

As detailed in Section \ref{sec:PM}, the PM method can produce an accurate approximation of the quasicrystals by solving the LP model in a higher dimensional reciprocal space while taking advantage of FFT to deal with the linear and nonlinear terms. However, determining a delicate phase diagram of the LP model is still expensive due to the wide range of parameters and the existence of multiple stable states. Although adaptive method exists \cite{Softwarejiang} that can generate the phase diagram without having to resolve the full parameter domain, it fails to produce the field variables $\widehat{\phi}$ and ${\phi}$ which are needed, e.g., in controlling the self-assembly of quasicrystals and a variety of other desired structures in practical experimental realizations\cite{provatas2011phase,chen2002phase,barkan2011stability,lifshitz2007soft}.

The proposed MCMS-RBM strives to learn the parameter dependence of the PM solution, vastly accelerate the generation of the phase diagram, and furnish surrogate and equally accurate field variables anywhere in the parameter domain. The MCMS-RBM has two prominent features in comparison to the standard RBM. First, it has multiple components with each providing reduction functionality for one branch of the problem induced by one part of the parameter domain. Second, without resorting to the parameter domain configurations (e.g. phase diagrams) \textit{a priori}, it enriches each component in a greedy fashion via a phase-transition guided exploration of the multiple states inherent to the problem. Specifically, it tests each stable state for every parameter value and retain all solutions that have not gone through any phase transitions. All these solutions, that are multiple for each parameter value, are adopted by the RBM components according to their convergent state.

We devote the rest of this section to the presentation of the two parts of the MCMS-RBM, namely its online and offline procedures, the adaptive algorithm that we adopt from \cite{Softwarejiang}, and its enhancement by the MCMS-RBM.

A key strategy of the RBM is the offline-online decomposition. During the offline procedure, five low-dimensional RB spaces of dimensions $N_{\rm S}$
\[
W_{N_{\rm S}}^{\rm S} \mbox{ for each state } \rm S \in \{{\rm QC, C6, LQ, T6, Lam}\}
\]
are generated by a greedy algorithm. These are called the five \textit{components} of the MCMS-RBM. During the online procedure for any given parameter value $\bmu$, the unknown RB coefficients in each component are solved through a reduced order model with an initial value given in the corresponding \textit{state}. Here, we first introduce the online procedure which will be repeatedly called during the offline construction phase to build $\{W_{N_{\rm S}}^{\rm S}\}$. To achieve online-efficiency, we resort to the EIM \cite{barrault2004empirical,grepl2007efficient}.

\subsection{Empirical interpolation method}
\label{sec:eim}
There are two nonlinear terms, one in the LP model $g(\phi(\bm{x},t);\bmu)=\alpha \phi^2-\phi^3$ and one in the free energy functional $ h(\phi(\bm{x},t);\bmu)=-\frac{\alpha}{3}\phi^3+\frac{1}{4}\phi^4$. Via a greedy algorithm that identifies function-specific ($\bmu$-independent) interpolation bases and corresponding interpolation points, the EIM approximates both functions by their interpolants
\begin{align*}
g(\phi(\bm{x},t);\bmu)& \approx g_{M_g}(\phi(\bm{x},t),t;\bmu)= \sum_{m=1}^{M_g} d_{M_g,m}(\bmu,t)g(\bm{x};\bmu_g^m),\\
h(\phi(\bm{x},t);\bmu)& \approx h_{M_h}(\phi(\bm{x},t),t;\bmu)= \sum_{m=1}^{M_h} \gamma_{M_h,m}(\bmu,t)h(\bm{x};\bmu_h^m).
\end{align*}
Here, $\{\bmu_g^m\}_{m=1}^{M_g}$ and $\{\bmu_h^m\}_{m=1}^{M_h}$ are two sets of parameter ensembles chosen by the greedy algorithm, while $g(\bm{x};\bmu_g^m)$ and $h(\bm{x};\bmu_h^m)$ are the corresponding functions that are orthonormal under point evaluations at the interpolation points. 

Specifically to our developed MCMS-RBM, we  construct one set of EIM expansions for each of the five states. Indeed, for each parameter initiated with each of the five states, the high-fidelity solution that has not gone through any phase transitions in the evolution will be adopted in the greedy procedure.

\subsection{Online procedure}

The MCMS-RBM approximates the high fidelity solution with a surrogate solution
\begin{equation}
\widehat{\phi}(\bmu,t)\approx\widehat{\phi}_{\text{rb}}(\bmu,t)=W_N c_N(\bmu,t),
\end{equation}
where $W_N$ is one of the $W_{N_{\rm S}}^{\rm S}$'s. For simplicity, we omit the sub- and sup-scripts $\rm S$ whenever there is no confusion. Moreover, $c_N$ is the RB coefficient to be solved for that component whose notational dependence on $\rm S$ is also omitted. 
The inverse Fourier transform of this surrogate solution can be directly derived as 
\begin{align*}
\phi (\bmu, t) \approx iW_Nc_N(\bmu, t),
\end{align*}
with $iW_N$ representing the inverse Fourier transform of basis space $W_N$. Therefore, the 
repeated transformations in real and reciprocal spaces during the iteration will also be calculated in low dimensional space and it is independent of the degrees of freedom of the full order model.
For simplicity, we denote the unknown RB coefficients $c_N(\bmu,t)$ at iteration $t$ by $c_t(\bmu)$.
Substituting these surrogate solutions and the EIM approximation into Eq. \eqref{eq:semi_implicit_pm}, and using the Galerkin projection method, one can derive the reduced order model
\begin{equation}
\begin{aligned}
W_N^T\left(1+\Delta t \cdot c\left(1-\boldsymbol{k}^2\right)^2\left(q^2-\boldsymbol{k}^2\right)^2 \right)W_Nc_{t+1}
=W_N^T(1+\Delta t \cdot \varepsilon) W_Nc_t+\Delta t \cdot W_N^T\mathscr{F}(V_M^gd_t^M),
\end{aligned}
\end{equation}
where $\bk^2 = \sum_{k=1}^d g_k^2$. Here, $V_M^g$ is the basis space constructed by the EIM for $g(\phi(\bm{x},t);\bmu)=\alpha \phi^2-\phi^3$, and $d_t^M$ is the coefficient of the EIM at every iteration.
We rewrite this equation as 
\begin{equation}
\begin{aligned}
A_1c_{t+1}
=A_2c_t + \Delta t \cdot \varepsilon A_2c_t + A_3d_t^M,
\end{aligned}
\label{eq:rb:online}
\end{equation}
where $A_1 = W_N^T\left(1+\Delta t \cdot c\left(1-\boldsymbol{k}^2\right)^2\left(q^2-\boldsymbol{k}^2\right)^2 \right)W_N \in \mathbb{R}^{N\times N}$, $A_2 =W_N^TW_N \in \mathbb{R}^{N\times N}$ and $A_3 = \Delta t W_N^T\mathscr{F}(V_M^g) \in \mathbb{R}^{N\times M}$. 
We note that they can all be pre-computed during the offline procedure via updates as each snapshot is identified. Further, they also depend on $\rm S$ but with notational dependence omitted.

The calculation of the free energy can also be accelerated since we have
\begin{equation}
    E_{\rm rb}= \frac{c}{2}c_t^TE_1c_t-\frac{\varepsilon}{2}c_t^TA_2c_t- \left(\gamma_t^L\right)^T V_L^h \bs{1},
\end{equation}
where $ E_1=iW_N^T\left(\left(1-\boldsymbol{k}^2\right)^2\left(q^2-\boldsymbol{k}^2\right)^2\right)^2iW_N\in \mathbb{R}^{N\times N} $, $V_L^h$ is the basis space of EIM for nonlinear terms of the free energy functional $ h(\phi(\bm{x},t);\bmu)=-\frac{\alpha}{3}\phi^3+\frac{1}{4}\phi^4$,  $\gamma_t^L$ is the coefficient of EIM at terminal time, and $\bs{1}$ is a column vector of $1$'s. The detailed online algorithm is provided in Algorithm \ref{alg:Algorithm:online}.
\begin{algorithm}[H]
  \caption{Online procedure of the MCMS-RBM}  
  \label{alg:Algorithm:online}  
  \begin{algorithmic}[1] 
   \vspace{0.5ex}
\State \textbf{Input:} Pre-computed reduced matrices $A_1$, $A_2$, $A_3$, and $E_1$ for a given $\rm S \in \{{\rm QC, C6, LQ, T6, Lam}\}$, tolerance $\delta$, parameter $\bmu$, temporal step $\Delta t$ and total time $T$
\State Initialize the Fourier coefficient $\widehat{\phi}_0$ for state $\rm S$, and derive the initial RB coefficient $c_0$. Set $\text{res}=1$, and $t=0$ 
\State \mbox{\textbf{While}} $\text{res}\geq \delta  \quad$ and  $t  \leq T$
\State \quad Calculate the EIM coefficient $d_t^M$ for $\bmu$
\State \quad Solve the RB coefficient $c_{t+\Delta t}$ by Eq. \eqref{eq:rb:online}, and calculate $\text{res} =||c_{t+\Delta t} -c_{t}||_2$
\State \quad $t =t +\Delta t$
\State \mbox{\textbf{End While}}
\State  Calculate the EIM coefficient $\gamma_t^L$ for $\bmu$
\State \textbf{Output:} Surrogate solution $\phi_{t, \text{rb}}=W_Nc_t$ , $\widehat{\phi}_{t, \text{rb}}=iW_Nc_t$, and the 
free energy functional
$E=\frac{c}{2}c_t^TE_1c_t-\frac{\varepsilon}{2}c_t^TA_2c_t-\left(\gamma_t^L\right)^T V_L^h \bs{1}$
  \end{algorithmic}  
  \label{alg:mcms_online}
\end{algorithm}

The complexity of this procedure is
$$O(N_{t, \mathbf{RB}}\cdot(N^3+NM^2)),$$
where $N_{t, \bm{RB}}$ is the iteration times of the RBM, which is almost the same as the iteration times of solving the full order model. However, $O(M^2+N^3)$ is much smaller than $O(\mathcal{N}\log \mathcal{N})$ because $N^3 \ll \mathcal{N}$, and $NM^2 \ll \mathcal{N}$. 
Moreover, the surrogate solution obtained by the RBM is essentially an approximation of the truth solution at each moment, so asymptotically the phase transition between the surrogate solution and the high-fidelity solution occurs almost simultaneously.

\subsection{Offline procedure}

 Now we present the greedy algorithm for constructing the RB spaces. The following error indicator for the LP model is adopted
\begin{align*}
\Delta_n(\bmu):=\frac{||r_{n}(\bmu)||_2}{\beta(\boldsymbol{\mu})},
\end{align*}
similar to the traditional residual-based error estimators \cite{rozza2008reduced,dede2010reduced,
grepl2005posteriori,huynh2013static}. Here
\begin{align*}
\beta(\bmu)=\sigma_{\min}\mathcal{A}_1,
\end{align*}
is the smallest singular value of matrix $\mathcal{A}_1$ with $\mathcal{A}_1 =c\left(1-\boldsymbol{k}^2\right)^2\left(q^2-\boldsymbol{k}^2\right)^2 \in \mathbb{R}^{\calN \times n}$ and $I$ being the identity matrix. The residual $r_n$ is defined as 
\begin{align*}
r_n(\bmu) = \mathcal{A}_1c_t -\varepsilon W_nc_t - \mathscr{F}^{-1}(V_M^g)d_t^M.
\end{align*}
Note that the computational complexity of the error indicator is made independent of $\calN$ via an offline-online decomposition. Indeed, one has that 
\begin{equation}
\begin{aligned}
||r_n(t;\bmu)||_2^2&=\left(\mathcal{A}_1c_t -\varepsilon W_nc_t - \mathscr{F}^{-1}(V_M^g)d_t^M\right)^T\left(\mathcal{A}_1c_t -\varepsilon W_nc_t - \mathscr{F}^{-1}(V_M^g)d_t^M \right)\\
&=c_{t+1}^TB_1c_t -2\varepsilon c_t^TB_2c_t -2c_t^TB_3d_t^M +\varepsilon^2 c_t^TB_4c_t +2\varepsilon c_t^TB_5d_t^M +(d_t^M)^T B_6d_t^M,
\end{aligned}
\label{eq:rb:errorestimator}
\end{equation}
where $B_1 = \calA_1^T\calA_1 \in \mathbb{R}^{n\times n}$, $B_2 =\calA_1^TW_n \in \mathbb{R}^{n\times n}$, $B_3 =\calA_1^T\mathscr{F}^{-1}(V_M^g) \in \mathbb{R}^{n\times M}$, $B_4 =W_n^TW_n\in \mathbb{R}^{n\times n}$, $B_5 =W_n^T\mathscr{F}^{-1}(V_M^g) \in \mathbb{R}^{n\times M}$, and $B_6 =(\mathscr{F}^{-1}(V_M^g))^T\mathscr{F}^{-1} \in \mathbb{R}^{M \times M}$ can be pre-computed and updated at each greedy loop. Afterwards, the computation of the error indicator only depends on the numbers of the EIM and RB basis.

We are now ready to describe the greedy algorithm for constructing the five RB spaces $W_{N_{\rm S}}^{\rm S}$ with $\rm S \in \{{\rm QC, C6, LQ, T6, Lam}\}$.  We first discretize the parameter domain $\mathcal{D}$ by a sufficiently fine training set $\Xi_{\text{train}}$. For any given $\rm S$, the guiding principle is that the RB space $W_{N_{\rm S}}^{\rm S}$ should contain snapshots having the structure corresponding to the component $\rm S$. Thus unlike the vanilla RBM, the first snapshot cannot be totally random. 
Indeed, this can be realized by performing Algorithm \ref{alg:Algorithm:PMsolver} for several parameters until the first $\bmu$ whose PSS \eqref{eq:mu_pss} contains a branch corresponding to the current component $\rm S$ is identified. 
Next, we call the online solver for $\rm S$ via Algorithm \ref{alg:mcms_online} and calculate the error estimator for each parameter in the training set. 
The \textit{temporary candidate} for the $(n+1)^{\rm th}$ $(n=1,2,\cdots,N-1)$ parameter is selected as the maximizer of the error estimators
$$\boldsymbol{\mu}_{n+1}^{\rm temp}=\mathop{\arg \max}\limits_{\boldsymbol{\mu} \in \Xi_{\text {train}}} \Delta_n(\boldsymbol{\mu}).$$
We then call Algorithm \ref{alg:Algorithm:PMsolver} to obtain the high-fidelity PSS $\widehat{\Phi}(\cdot; \bmu_{n+1}^{\rm temp})$. If this PSS contains a branch corresponding to $\rm S$, we enrich $W_{N_{\rm S}}^{\rm S}$ with this branch\footnote{In practice, we incorporate a Gram-Schmidt procedure for numerical robustness.}. If $\widehat{\Phi}(\cdot; \bmu_{n+1}^{\rm temp})$ contains no $\rm S$-specific component, we prune $\boldsymbol{\mu}_{n+1}^{\rm temp}$ and go to the next maximizer of the error indicators. 
This is repeated until we find a $\boldsymbol{\mu}_{n+1}^{\rm temp}$ whose PSS has a $\rm S$-specific component. Finally, we set $\bmu_{n+1}$ to be this $\bmu_{n+1}^{\rm temp}$.
The detailed greedy algorithm for the construction of the RB spaces $\{W_{N_{\rm S}}^{\rm S}\}$ is described in 
Algorithm \ref{alg:Algorithm:offline}.
\begin{algorithm}[H]
  \caption{Offline procedure of MCMS-RBM}  
  \label{alg:Algorithm:offline}  
  \begin{algorithmic}[1]  
  \vspace{0.5ex}
\State \textbf{Input:} training set $\Xi_{\text{train}}$, EIM basis space $V_M^g$, collocation set $X_M$, component $\rm S$ to be built and tolerance $\delta$ 
\State Select $\bmu_1$ from $\Xi_{\text {train }}$ so that PSS $\widehat{\Phi}(\cdot; \bmu_{1})$ has an $\rm S$-state denoted by $\widehat{\phi}(\bmu_1)$. Initialize $S_1={\bmu_1}$, 
$W_1=\{\widehat{\phi}(\bmu_1)/||\widehat{\phi}(\bmu_1)||_2\}$, $iW_{1}= \mathscr{F}^{-1}(\widehat{\phi})/||\mathscr{F}^{-1}(\widehat{\phi})||_2$, and $n=1$
    \For {$n=2,3,\cdots,N-1$}
\State Solve the RB coefficient in the $\rm S$-component of MCMS-RBM, $c_n(\bmu)$, by Eq. \eqref{eq:rb:online} 

for each {$\bmu\in \Xi_{\text{train}}$}, and calculate the error estimator through Eq. \eqref{eq:rb:errorestimator}
    \State Sort $\Delta _n(\bmu)$ and select the first $\bmu$ whose PSS $\widehat{\Phi}(\cdot; \bmu_{1})$ has an $\rm S$-state denoted 
    
    by $\widehat{\phi}_{n+1}$, and set the corresponding $\bmu$ as the $(n+1)^{\rm th}$ parameter $\bmu_{n+1}$
\State Update $S_{n+1}=\{S_n,\bmu_{n+1}\}$, $W_{n+1}=\text{GS}(W_n,\widehat{\phi}_{n+1})$, $\phi_{n+1}=\mathscr{F}^{-1}(\widehat{\phi}_{n+1})$,
$\eta_{n+1}=$

$\phi_{n+1}/||\phi_{n+1}||_2,\quad$ $iW_{n+1}=\{iW_n,\eta_{n+1}\}$ and the reduced matrices $\{A_i\}_{i=1}^3,$ $ \{B_i\}_{i=1}^6$

 $\mbox{ and } E_1$
 \EndFor
 \State \textbf{Output:} Reduced basis space $W_n$, the inverse RB space $iW_n$, and the reduced matrices $\{A_i\}_{i=1}^3, \{B_i\}_{i=1}^6\, \mbox{ and } E_1$
  \end{algorithmic}  
\end{algorithm}

\subsection{Adaptive phase diagram generation by the MCMS-RBM}
\label{sec:adaptivephase}

One advantage of building an efficient surrogate solver such as the developed MCMS-RBM is that it makes feasible the generation of phase diagrams via an exhaustive sampling of the parameter domain. 
We propose to further leverage the efficiency of the MCMS-RBM by adopting the adaptive strategy \cite{Softwarejiang} when querying the parameter domain. 
The idea is to first sample a coarse cartesian grid to determine the phase of each point, and then for each point one checks the eight nearest neighbors. If any of the neighbors were deemed in a different phase, one regards the current point as a boundary point. In this situation, one appends the discrete grid with the middle point of the two points with different phases. We remark that, as the adaptive algorithm proceeds, the grid will become more unstructured making the ``eight nearest neighbors'' not as easily identifiable as the initial structured grid. In this case, we simply sort the neighbors by its Euclidean distances to the base point.

\section{Numerical results}

\label{sec:numerics}

In this section, we test the proposed MCMS-RBM on the two-dimensional quasiperiodic LP model parameterized by the reduced temperature $\varepsilon$ and the phenomenological parameter $\alpha$ delineating the level of asymmetry. Furthermore, we highlight its efficiency and accuracy by adopting the adaptive phase diagram generation algorithm to produce a phase diagram that is as accurate as the state of the art. 

\subsection{Setup and notations}
The parameter domain is set to be $\mathcal{D}=[-0.0125,0.0515]\times [0,1]$. The training and testing sets are $\mathcal{D}$'s two disjoint uniform cartesian discretizations, 
$$
\begin{aligned}
& \Xi_{\text{train}}=(-0.0125:0.002:0.0515)\times(0:0.05:1),\\
& \Xi_{\text{test}}=(-0.01:0.015:0.05)\times(0.025:0.1:0.925).\\
\end{aligned}
$$
For the high-fidelity solver of the LP model, we set the degree of freedom of the Fourier spectral method in each direction as $N_{\mathbf{H}}=32$. 
We denote the worst-case relative errors of the nonlinear terms with $m$- dimensional space $V_M^g$ by $\phi_{\text{err}}^g(m)$, solutions and free energy functionals with $n$- dimensional space $W_n$ by $\phi_{\text{err}}(n)$ and $E_{\text{err}}(n)$, respectively.
\begin{align*}
&\phi_{\text{err}}^g(m)=\max _{\boldsymbol{\mu} \in \Xi_{\text {train }}}\left\{\frac{\left\|g(\bmu)- g_{M_g}(\bmu)\right\|_2}{\|g(\bmu)\|_2}\right\}, \\
&\phi_{\text{err}}(n)=\max _{\boldsymbol{\mu} \in \Xi_{\text {test }}}\left\{\frac{\left\|\widehat{\phi}(\bmu)-\widehat{\phi}_{n,rb}(\boldsymbol{\mu})\right\|_2}{\|\widehat{\phi}(\bmu)\|_2}\right\}, \\
&E_{\text{err}}(n)=\max _{\boldsymbol{\mu} \in \Xi_{\text {test }}}\left\{\frac{\left|\mathcal{F}(\widehat{\phi}(\boldsymbol{\mu}))-\mathcal{F}(\widehat{\phi}_{n,rb}(\boldsymbol{\mu}))\right|}{|\mathcal{F}(\widehat{\phi}(\bmu))|}\right\}.
\end{align*}
Finally, we denote by $\Delta_{RB}(n)$ the worst-case error indicators when the RB spaces are of $n$-dimensional,
$$\Delta_{RB}(n)=\max_{\bmu\in \Xi_{\text{train}}}\Delta_n(\bmu).$$

\subsection{MCMS-RBM results}
We are now ready to present the numerical results of the MCMS-RBM applied to the parameterized two-dimensional quasiperiodic LP model.

\begin{figure}[t!]
\centering
\includegraphics[width=0.8\textwidth]{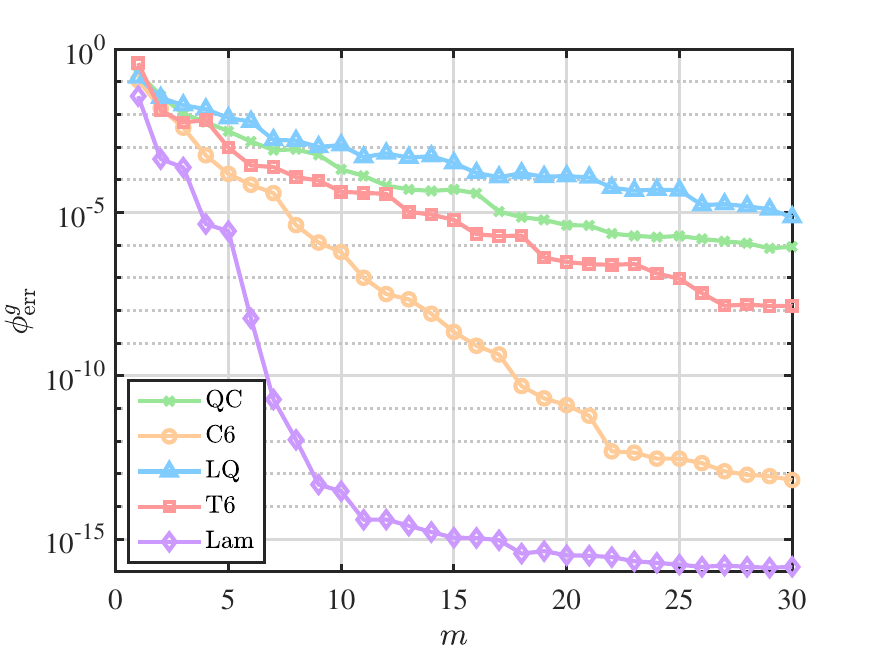}
\caption{Relative EIM errors}
\label{fig:EIM}
\end{figure}

\begin{figure}
  \begin{subfigure}[b]{0.32\textwidth}
    \centering    \includegraphics[width=\textwidth]{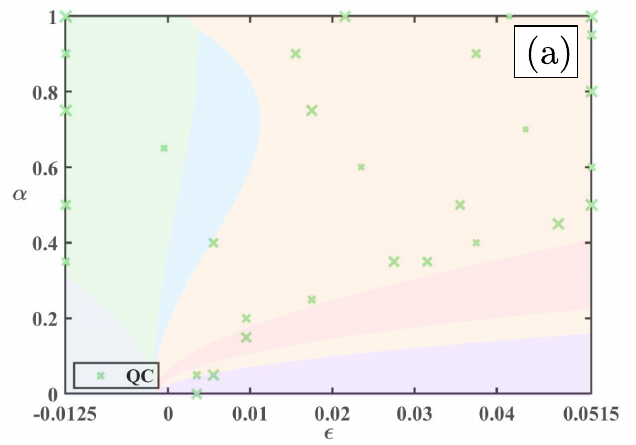}
  \end{subfigure}
  \hfill
  \begin{subfigure}[b]{0.32\textwidth}
    \centering    \includegraphics[width=\textwidth]{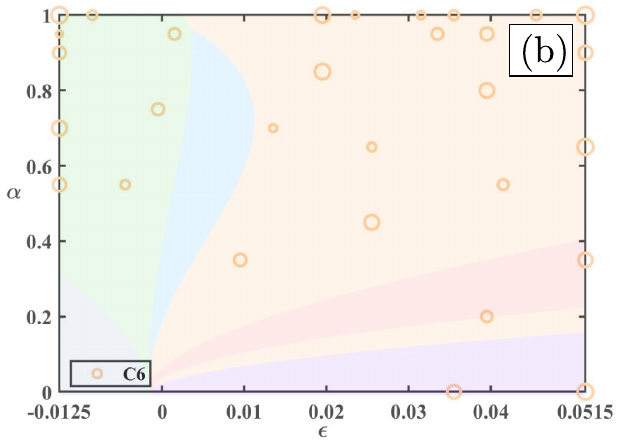}
  \end{subfigure}
  \hfill
  \begin{subfigure}[b]{0.32\textwidth}
    \centering    \includegraphics[width=\textwidth]{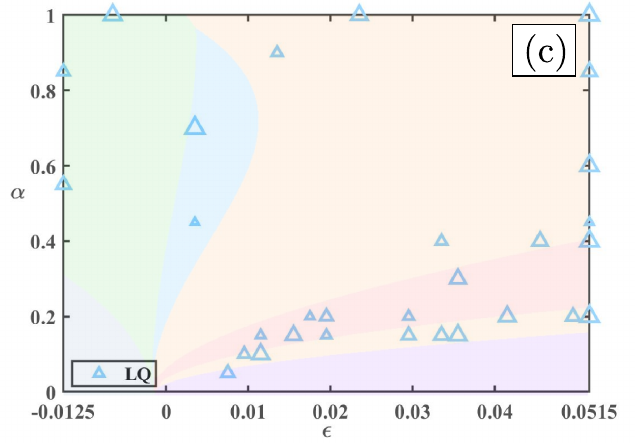}
  \end{subfigure}

  \begin{subfigure}[b]{0.32\textwidth}
    \centering    \includegraphics[width=\textwidth]{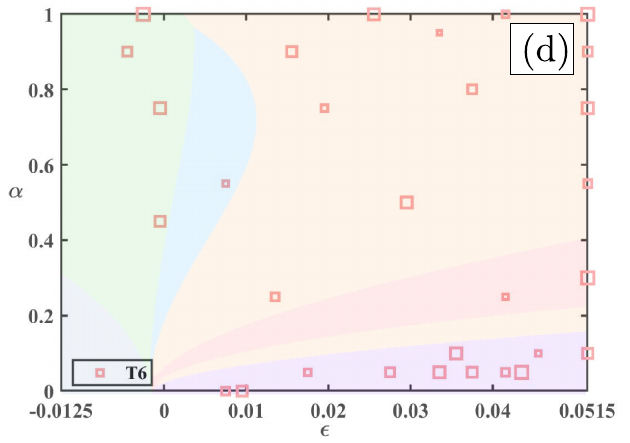}
  \end{subfigure}
  \hfill
  \begin{subfigure}[b]{0.32\textwidth}
    \centering    \includegraphics[width=\textwidth]{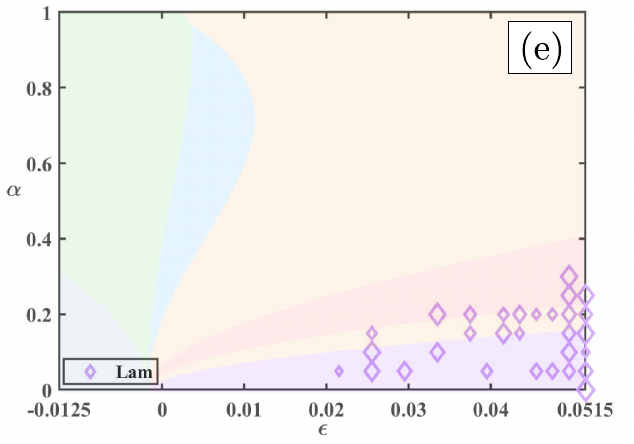}
  \end{subfigure}
  \hfill
  \begin{subfigure}[b]{0.32\textwidth}
    \centering    \includegraphics[width=\textwidth]{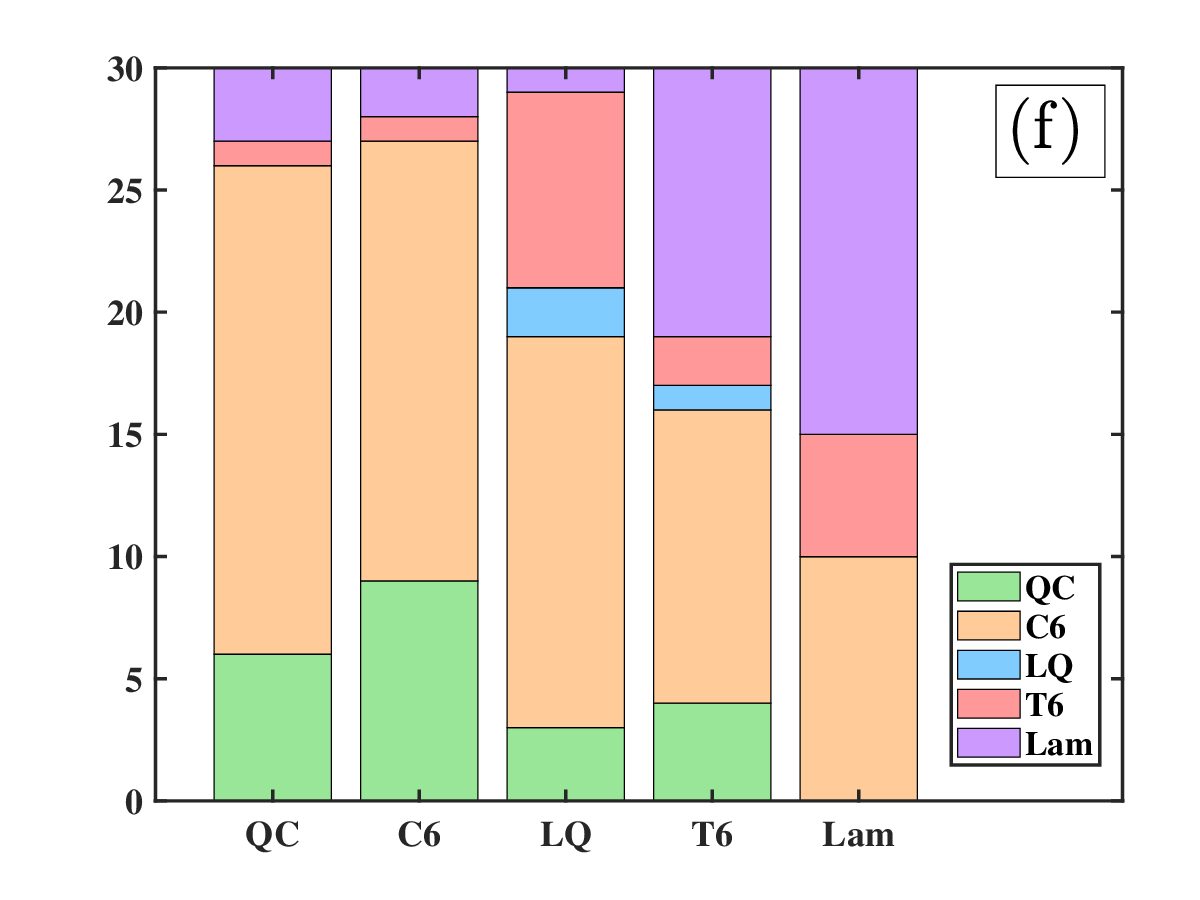}
  \end{subfigure}
\caption{Distributions of the parameters selected by EIM corresponding to phases QC, C6, LQ, T6 and Lam (a-e) , and the histogram of each set of basis (f). }
\label{fig:para:EIM}
\end{figure}

{\bf EIM results} ---
The relative error curves of the nonlinear terms in the LP model and the distribution of the first $M(=30)$ selected parameters for each of the five MCMS-RBM components are showed in Fig. \ref{fig:EIM} and Fig. \ref{fig:para:EIM}. As expected, the error curves exhibit exponential convergence as the number of the EIM basis increases. It is worth noting that all five states can effectively limit the error to within $10^{-4}$ using just up to $30$ basis functions. This allows for significant speedup for the reduced model. As to the distribution of the chosen parameter values, the selected parameters predominantly lie within the corresponding state of their phase diagrams for relatively simple structures like Lam 
(Fig. \ref{fig:para:EIM}(e)). However, for other complex structures such as QC
and LQ (Fig. \ref{fig:para:EIM}(a, c)), some parameters are chosen from other states and with a more uniform distribution with clusters toward the boundary. In Fig. \ref{fig:para:EIM}(f), the histogram is displayed to show the number of parameters distributed in different states. These results underscore \textit{the many-to-many feature between components and states} of the developed MCMS-RBM.

\begin{figure}[H]
\includegraphics[scale=0.48]{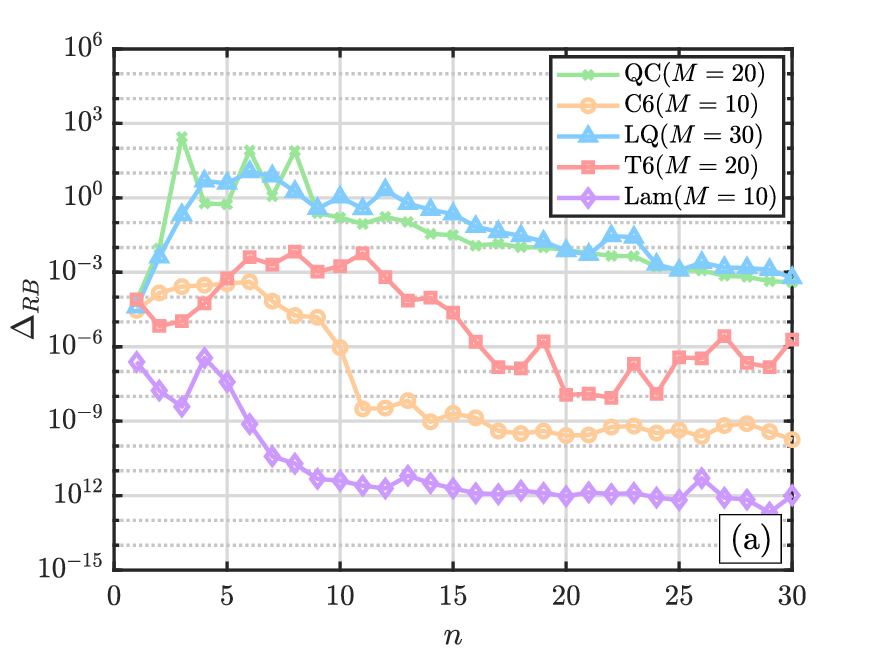}
\includegraphics[scale=0.48]{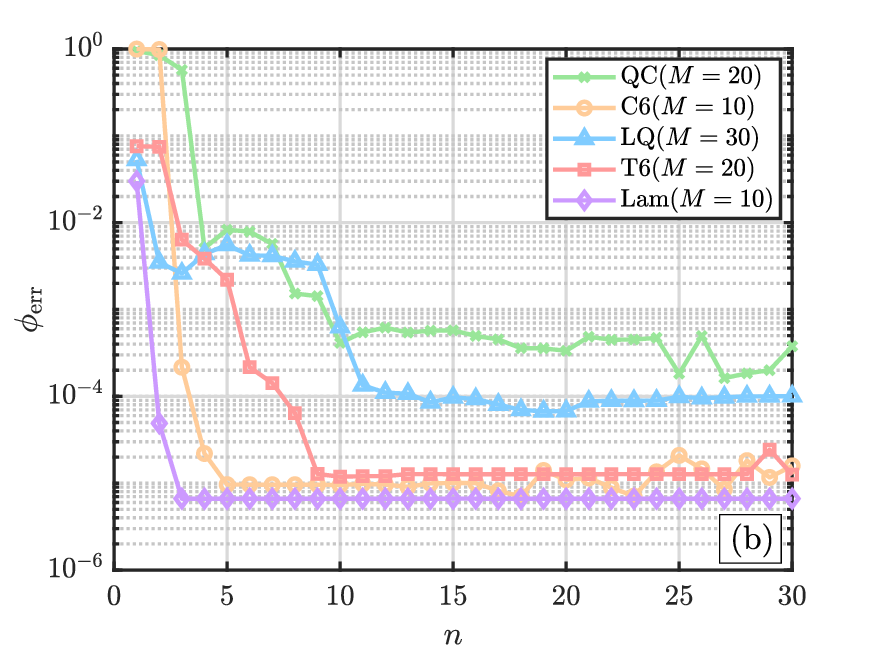}

\includegraphics[scale=0.48]{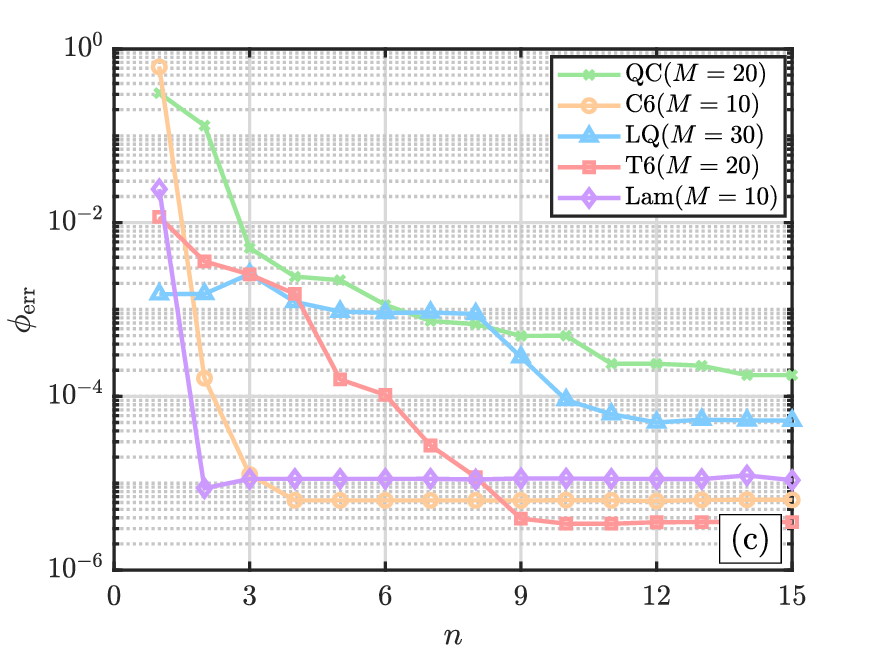}
\includegraphics[scale=0.48]{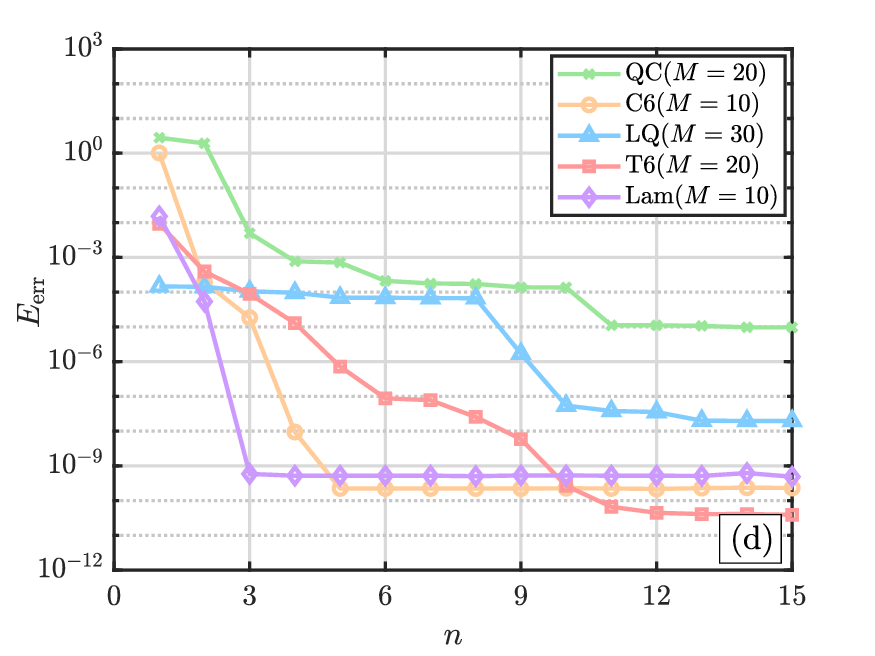}
\caption{Training (a-b) and testing (c-d) results of MCMS-RBM. Top: Error indicators (a) and relative error curves during training  (b). Bottom: Relative testing error curves of the solutions (c) and free energy functional (d).}
\label{RBM_results}
\end{figure}

\begin{figure}
  \centering
  \begin{subfigure}[b]{0.32\textwidth}
    \centering    \includegraphics[width=\textwidth]{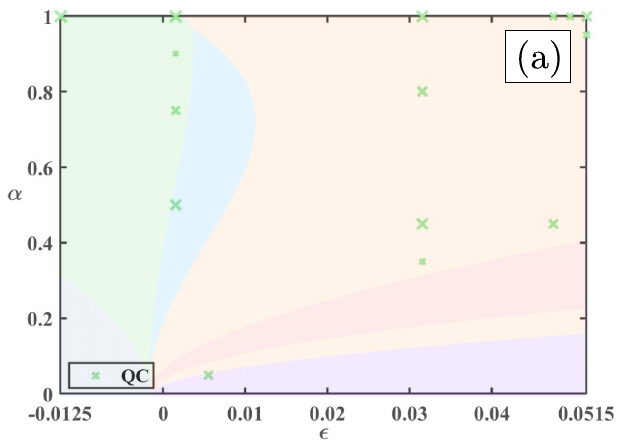}
  \end{subfigure}
  \hfill
  \begin{subfigure}[b]{0.32\textwidth}
    \centering    \includegraphics[width=\textwidth]{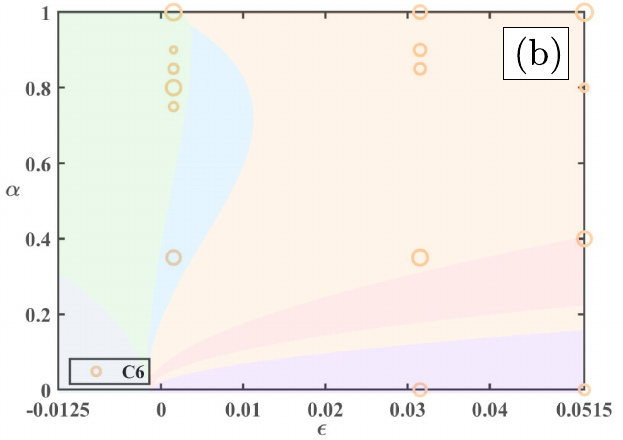}
  \end{subfigure}
  \hfill
  \begin{subfigure}[b]{0.32\textwidth}
    \centering    \includegraphics[width=\textwidth]{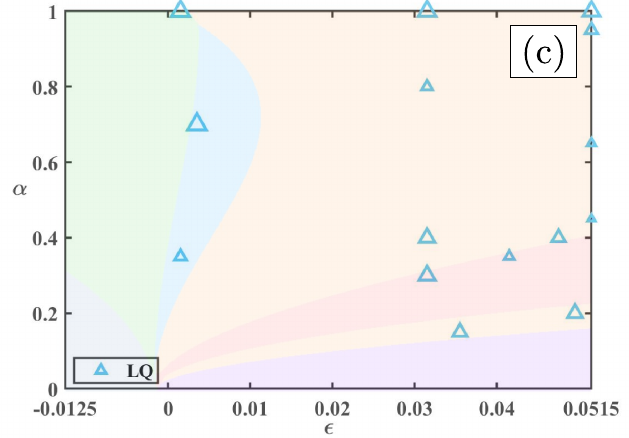}
  \end{subfigure}

  \begin{subfigure}[b]{0.32\textwidth}
    \centering    \includegraphics[width=\textwidth]{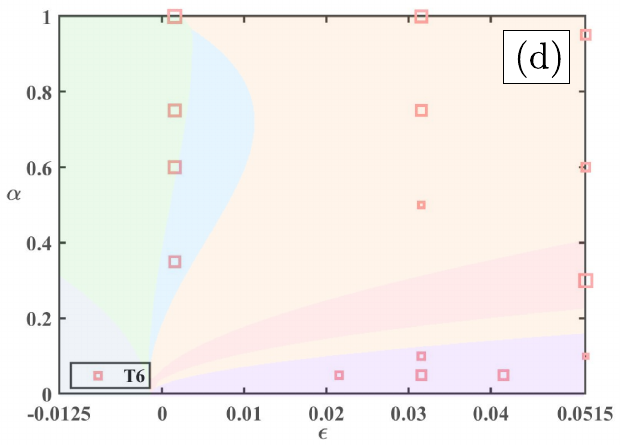}
  \end{subfigure}
  \hfill
  \begin{subfigure}[b]{0.32\textwidth}
    \centering    \includegraphics[width=\textwidth]{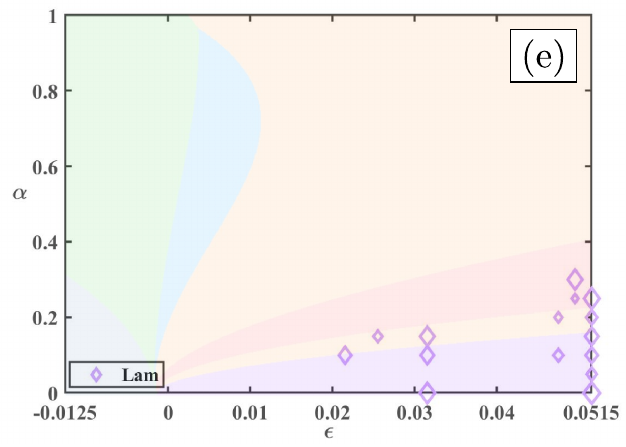}
  \end{subfigure}
  \hfill
  \begin{subfigure}[b]{0.32\textwidth}
    \centering    \includegraphics[width=\textwidth]{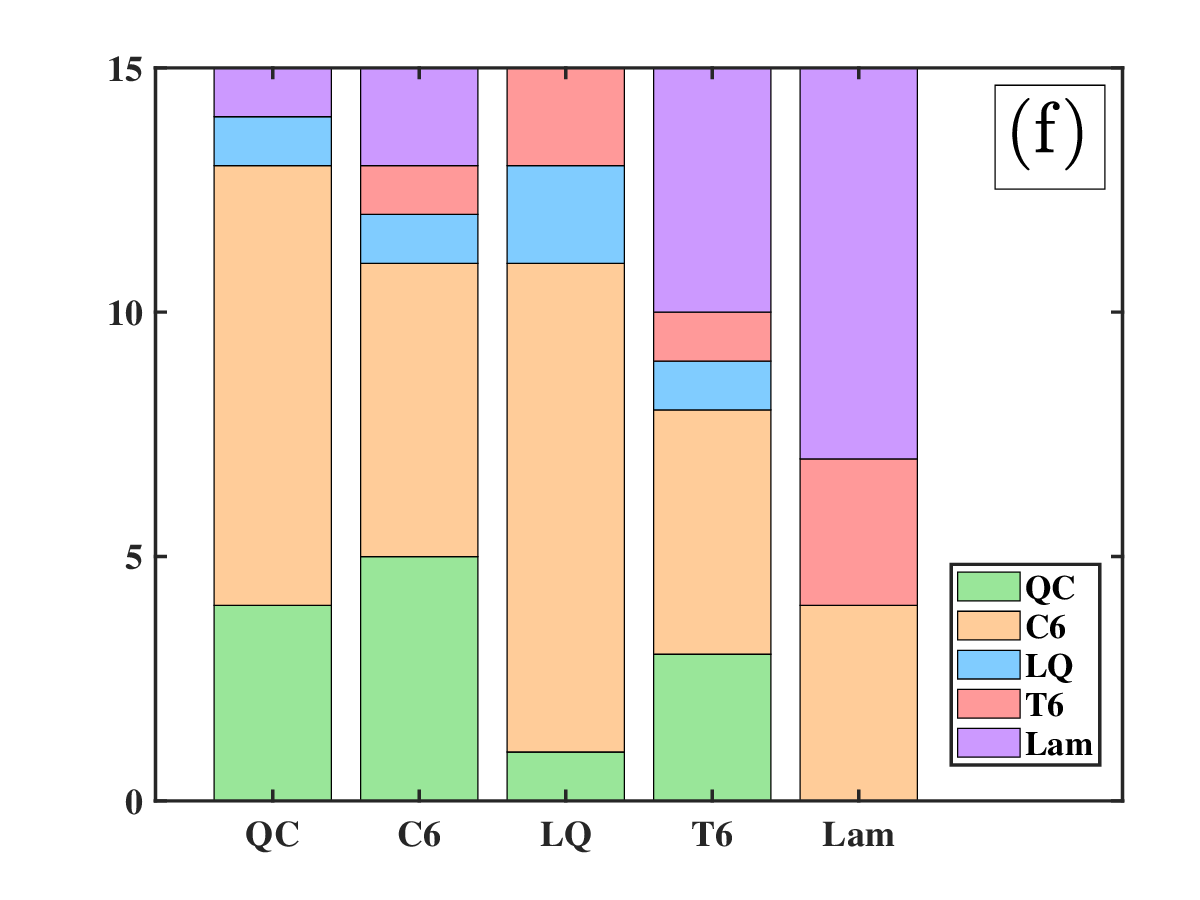}
  \end{subfigure}
\caption{Distributions of the parameters selected by MCMS-RBM corresponding to phases QC, C6, LQ, T6 and Lam (a-e) , and the histogram of each set of basis (f). }
\label{fig:para:offline}
\end{figure}

 {\bf RBM results} --- In Fig. \ref{RBM_results} (a, b), we present the error indicators and relative errors during the offline training stage. Initially, when the number of RB basis is insufficient, the corresponding error estimator increases. However, all the curves demonstrate exponential convergence rates when the reduced basis spaces are sufficiently enriched. The efficiency of the error indicator and the accuracy of the MCMS-RBM are underscored by the observed stable exponential convergence that the relative error curves exhibit right from the beginning. Furthermore, the relative errors of the solution with only $15$ basis can reach $10^{-3}$. The testing errors, in both the solution and the functional, are shown in Fig. \ref{RBM_results} (c, d). They, too, decrease exponentially. 
The distribution of the selected parameters of the RBM is shown in Fig. \ref{fig:para:offline}. The pattern resembles that of the EIM process.

Finally, we select five different parameters with one from each phase and apply our MCMS-RBM online solver to the corresponding LP model. We present in Table \ref{table:time} the EIM/RB dimensions, the wall clock times for the full problem and the reduced problem, and the MCMS-RBM errors in the order parameter $\phi$ and the energy functional ${\mathcal F}$. The MCMS-RBM consistently achieves an acceleration three orders of magnitude while both relative errors are at levels of $10^{-5}$ and $10^{-9}$. 
\begin{table}[!htb]
	\centering
\setlength{\tabcolsep}{5pt}
\renewcommand{\arraystretch}{1.5}
			\begin{tabular}{cccccc}
		\hline
		$\bmu$	&$(5\times 10^{-6},{\sqrt{2}}/{2})$ &$(0.05,1)$ &$(0.005,0.6)$ &$( 0.05,0.3)$&$(0.05,0.1)$
 \\ \hline
Phase&~~QC~~&C6 &  LQ &T6 & Lam \\ \hline
$(M,N)$& (20,15)&(10,5)&(30,15)&(20,10)&(10,5)\\ \hline
FOM time	&1.06e+02 &1.46e+01  & 7.95e+01 &3.85e+01&5.02e+01\\  \hline
RBM time	&2.98e-02    &1.31e-02 &4.14e-02 &2.28e-02&1.79e-02\\ \hline
$\phi_{\text{err}}$	&4.67e-06    &3.02e-06 &3.53e-05&3.17e-05&2.75e-05\\ \hline
$E_{\text{err}}$	&2.93e-10    &4.52e-11 &7.70e-09&1.83e-09&3.02e-09\\ \hline
	\end{tabular}
	\caption{Online computation time (in seconds) and relative error for the LP model at different parameters with $\Delta t=0.1$ and $N_{\bm{H}}=32$.}
	\label{table:time}
\end{table}

\subsection{Phase diagram generation}

\begin{figure}[H]
\includegraphics[scale=0.46]{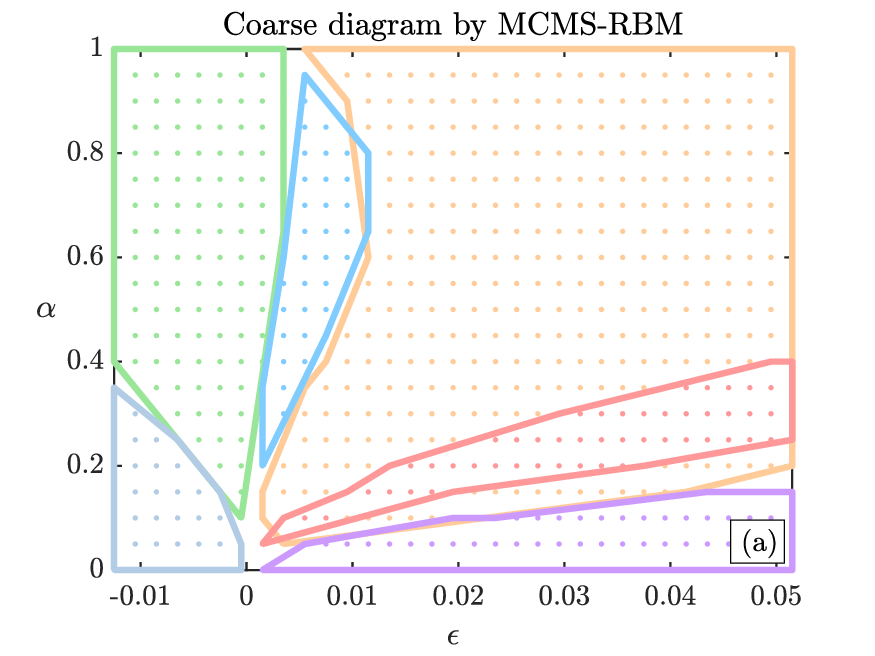}
\includegraphics[scale=0.46]{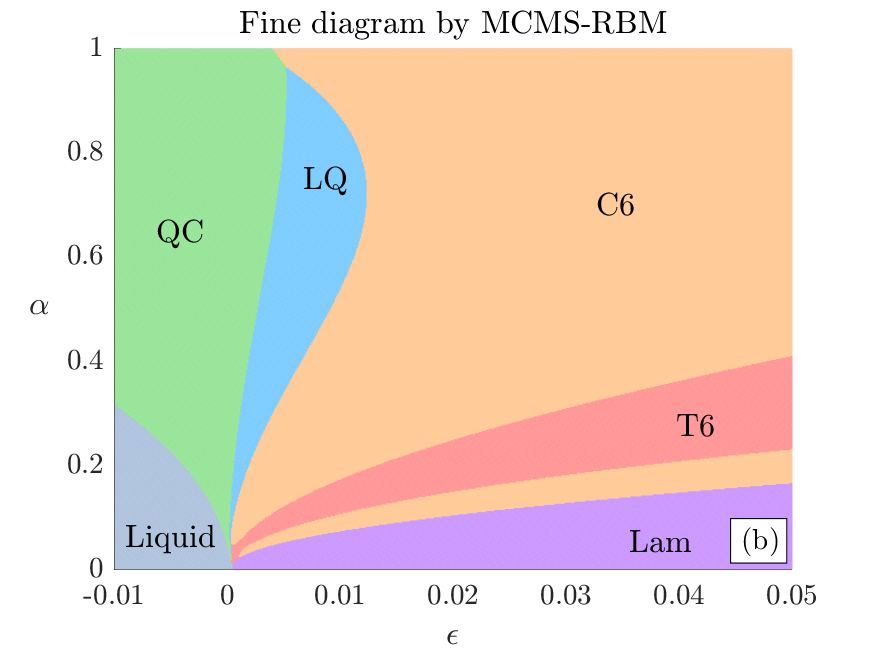}

\includegraphics[scale=0.31]{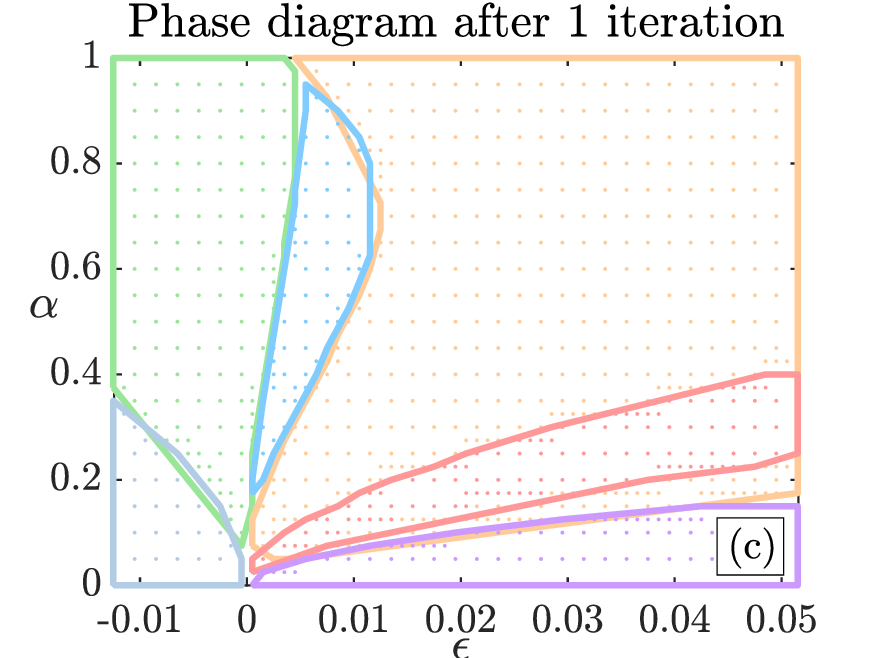}
\includegraphics[scale=0.31]{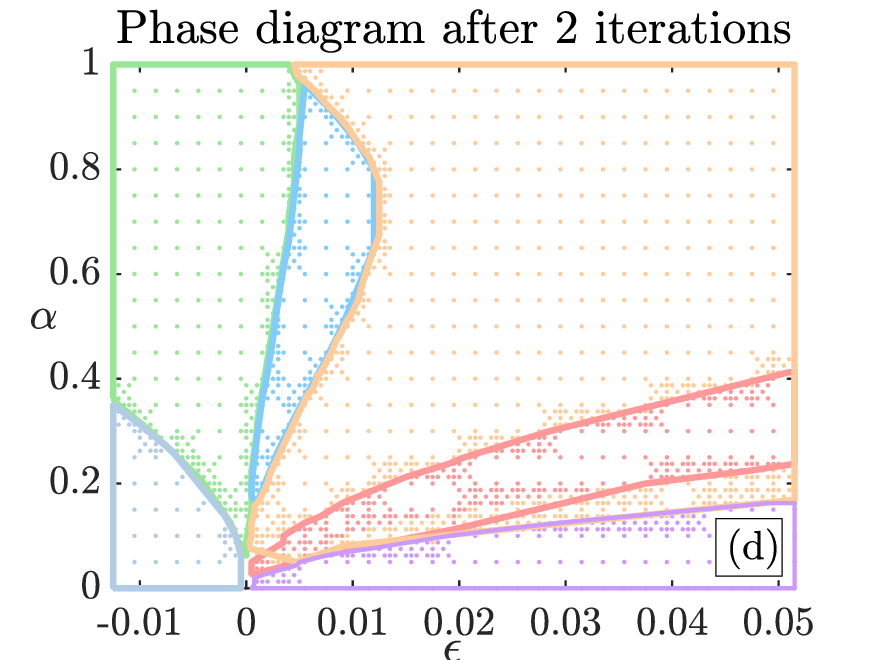}
\includegraphics[scale=0.31]{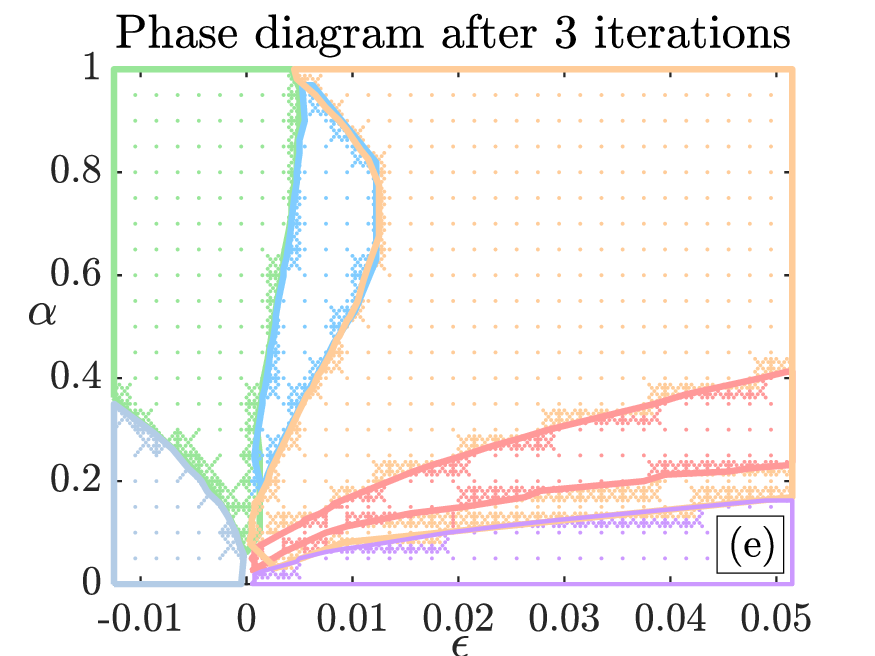}
\caption{The iterative phase diagrams (a,c-e) generated by the adaptive MCMS-enhanced algorithm. Shown on top right (b) is a phase diagram with high resolution. }
\label{fig:adaptive}
\end{figure}

We now apply the adaptive phase generation algorithm of Section \ref{sec:adaptivephase} to our parametric LP model with the EIM and RB dimensions given as in Table \ref{table:time}. The initial phase diagram on a uniform coarse grid  is generated by repeatedly invoking the online solver of the MCMS-RBM, see Fig. \ref{fig:adaptive} (a). Then the adaptive refinement is performed along the automatically detected boundaries of adjacent phases. The bottom row of Fig. \ref{fig:adaptive} contains the results of three consecutive refinements, with the third one capturing the delicate boundaries of the phase diagram quite well. 
In comparison, we query MCMS-RBM on a highly refined discretization of the parameter domain, a $601\times 501$ uniform grid. See Fig. \ref{fig:adaptive} (b) for the resulting phase diagram. It is clear that the third iteration of the adaptive algorithm agrees with this fine phase diagram which is in turn consistent with that in \cite{yin2021transition}.

As evidence of the efficiency of the adaptive algorithm, we list the wall clock time (in seconds) of each iteration in Table \ref{tab:adaptime}. It is clear that the generation of the coarse diagram  only takes $120.31$ seconds, and this value will increase when one performs one more iteration, as shown in the 7th column of the table. The total time of the generation of Fig. \ref{fig:adaptive}(e) is 276.65, which is derived by summarizing all the values of the 7th column.  Indeed, the $601\times 501$ fine phase diagram takes about $16$ hours in a serial environment by MCMS-RBM. A simple scaling with the acceleration rate of Table \ref{table:time} indicates that this phase diagram would have taken over $40$ months to generate if we were to call the FOM solver repeatedly. In a word, the computational cost of the phase diagram is significantly reduced by the MCMS-RBM and it can be further reduced by the adaptive refine boundary algorithm. 
\begin{table}[!htb]
\centering
\setlength{\tabcolsep}{5pt} 
\renewcommand{\arraystretch}{1.5} 
\begin{tabular}{cccccccc}
\hline
Phase&QC&C6 &LQ &T6 &Lam &Total \\ \hline
coarse &33.94&11.76&41.77&20.33&12.50
&120.31\\ \hline
1st iteration	&3.88 &3.97  &4.88 &4.11&1.32&18.17\\  \hline
2nd iteration	&11.61&10.51 &13.90 &11.00&2.18&49.20\\ \hline
3rd iteration	&26.46&17.27 &23.47 & 18.18&3.59&88.97\\ 
 \hline
\end{tabular}
\caption{Online computation time (in seconds) for adaptive refine boundary algorithm. }
\label{tab:adaptime}
\end{table}

\section{Conclusion}
\label{sec:conclusion}

This paper proposes a multi-component multi-state reduced basis method (MCMS-RBM) for the parametrized quasiperiodic LP model with two length scales. Featuring multiple components with each providing a reduced order model for one branch of the problem induced by one part of the parameter domain, the MCMS-RBM serves as a generic framework for reduced order modeling of parametric problems whose solution has multiple states across the parameter domain. 

Via a greedy algorithm that identifies the representative parameter values and a phase transition indicator, the method searches for the (potentially multiple) phase-steady solutions for each parameter value which are then used to enrich the corresponding components of the MCMS-RBM. 
Numerical experiments corroborate that the method can provide surrogate and equally accurate field variables, with speedup of three orders of magnitude, anywhere in the parameter domain. It can also accelerate the generation of a delicate phase diagram to a matter of minutes.

\bibliographystyle{abbrv}
\bibliography{MCMS_LP.bbl}
\end{document}